\documentclass{amsart}
\usepackage[initials]{amsrefs} 

\theoremstyle{plain}
  \newtheorem{theorem}    {Theorem}[section]
  \newtheorem{corollary}  [theorem]{Corollary}
\numberwithin{equation}   {section}

\newcommand{\N}{{\bf N}}
\newcommand{\Q}{{\bf Q}}
\newcommand{\R}{{\bf R}}
\newcommand{\B}{{\bf B}}
\newcommand{\Oo}{{\bf O}}
\newcommand{\Tmin} {T_{\text{min}}}
\DeclareMathOperator{\codim}{codim}
\DeclareMathOperator{\Vol}  {Vol}

\makeatletter
\@namedef{subjclassname@2010}{%
  \textup{2010} Mathematics Subject Classification}
\makeatother

\begin{document}

\title
           [Zariski decomposition  and           Morse inequalities]
{Divisorial Zariski decomposition\\and algebraic Morse inequalities}
\author
[S.      Trapani]
{Stefano Trapani}
\address
{Dipartimento di Matematica\\
 Universit\`a di Roma ``Tor~Vergata''\\
 Viale della Ricerca Scientifica\\00133~Roma\\Italy}
\email{trapani@mat.uniroma2.it}
\begin{abstract}
In this note we use the divisorial Zariski decomposition to give a more precise version of the algebraic Morse inequalities.
\end{abstract}
\maketitle


\section{Introduction}
In~\cite{Dem-5} J.-P. Demailly proved the so-called holomorphic Morse inequalities. They relate the cohomology groups of a line bundle on a compact complex manifold to the distribution of the eigenvalues of the curvature of a smooth Hermitian metric on the given line bundle. In~\cite{Tra} we considered the difference $L-F$ of a nef line bundle $L$ and a pseudoeffective line bundle $F$ on a projective manifold in order to derive a more intrinsic version of some of his inequalities. Our statement did not make use of a specific smooth Hermitian metric but it used intersection numbers of the bundles. The same result, in the case both bundles are nef, was proved and used by Y.~T. Siu to derive his effective Matsusaka big theorem \cite{Siu}. A more precise version of this theorems was later obtained in the case $F$ nef by Demailly in~\cite{Dem-1}. Even though in~\cite{Tra} we considered the case $F$ pseudoeffective and not necessarily nef, if $F$ is not nef a singular Hermitian metric on $F$ comes into the picture, and this makes the results not really intrinsic. The purpose of this paper is to show that the concept of non-nef locus and the divisorial Zariski decomposition obtained in~\cite{Bou} seem to give the correct framework for the algebraic Morse inequalities.

\section{Preliminaries}
Let $X$ be a compact projective manifold of complex dimension $n$. A divisor $F$ on $X$ is said to be \textit{pseudoeffective} if $c_1(F)$ is in the closure in the real Neron-Severi space of the cone of classes of effective $\Q$ divisors. Alternatively, $F$ is pseudoeffective if and only if the line bundle associated to $F$ carries a singular Hermitian metric of non-negative curvature current. More generally a cohomology class in $H^2(X,\R)$ is said to be pseudoeffective if it is represented by a non-negative closed $(1,1)$-current. Let $T$ be a closed non-negative $(1,1)$-current on $X$, and denote by $\nu_x(T)$ the Lelong number of $T$ at $x$. Siu \cite{Siu-2} proved that for every positive number $c$ the set $E_c(T)=\{x\in X:\nu_x(T)\geq c\}$ is a proper closed analytic subset of $X$. Let $b_p(X)=\inf\{c>0:\codim E_c(T)\geq p\}$. For the necessary background on singular metrics, currents, Lelong numbers, currents of minimal singularities and so on, see for example \cites{Bou,Dem-1,Dem-2,Dem-3}.

The \textit{stable base locus} of $F$ is the intersection of base loci of tensor powers of $F$, that is $\B(F)=\bigcap_m Bs(mF)$. Let $F$ be a $\Q$ divisor and let $S=\{m\in\N:\text{$mF$ is a divisor}\}$. The \textit{volume} of $F$ is by definition given by
\begin{equation*}
\Vol(F)=\limsup_{\substack{m\to\infty\\m\in S}} \frac{n!\,h^0(X,mF)}{m^n}.
\end{equation*}
The volume can be naturally extended to a non-negative continuous function on the set of real pseudoeffective cohomology classes. Recall that the interior of the set of pseudoeffective classes is the set of big classes, i.e., the set of classes represented by a current larger then a K\"{a}hler form. Moreover a pseudoeffective class of strictly positive volume is big. Therefore the extension to $0$ of the volume function is continuous on all of $H^2(X,\R)$ \cite{Bou-2}.

Let $A$ be an ample line bundle on $X$. The \textit{non-nef locus} of $F$ is $E_{n,n}(F)=\bigcup_k \B(kF+A)$, and this set in independent of the choice of the ample bundle $A$. Then $F$ is \textit{numerically effective} (\textit{nef}) if and only if $E_{n,n}(F)=\emptyset$. Alternatively if $F$ is a pseudoeffective line bundle and $\epsilon$ a positive rational number, let $\Tmin(\epsilon)$ be a current of minimal singularity in the family of currents $T$ representing $c_1(F)$ such that $T\geq-\epsilon A$. Then $\nu_x(\Tmin(\epsilon))$ is decreasing in $\epsilon$. Let $\nu^*_x(F)=\sup_\epsilon \nu_x(\Tmin(\epsilon))$. Then $E_{n,n}(F)=\{x\in X:\nu^*_x(F)\neq 0\}$, see~\cite{Bou}. A similar notion can be given for the case of pseudoeffective cohomology classes. (The non-nef locus is also called the \textit{diminished base locus} and is denoted by $\B_- (F)$, see for example \cite{H-McK}).

The codimension of $E_{n,n}$ is defined by $\codim E_{n,n}=\min_k \codim\B(kF+A)$. So $F$ is nef if and only if $\codim E_{n,n}(F)=n+1$. Since the non-nef locus cannot be of dimension $0$ then $\codim E_{n,n}(F)\neq n$. Moreover for every pseudoeffective divisor $F$ the codimension of $E_{n,n}(F)$ is at least $1$, and it is at least $2$ if and only if $F$ is \textit{modified nef} in the language of~\cite{Bou}.

Let $F$ be a pseudoeffective line bundle. In~\cite{Bou} S. Boucksom proves the following facts. There exists a finite number $D_1,\dotsc,D_N$ of codimension-one components of the non-nef locus of $F$. Moreover if we let $\nu_j=\min_{x\in D_j} \nu^*_x(F)$, then $\nu_j>0$ for $1\leq j\leq N$ and the non-nef locus of the real divisor $Z(F)=F-\sum_j \nu_j D_j$ has codimension at least $2$. See also \cite{Nak}. The decomposition $F=Z(F)+\sum_j \nu_j D_j$ is called the \textit{divisorial Zariski decomposition} of $F$. Note that Boucksom's theorem is true in the more general context of real pseudoeffective classes.

Let $TX$ be the tangent bundle of $X$, and let $\pi:P(T^* X)$ be the projectivized cotangent bundle with its natural projection map. Let $\Oo_{TX} (1)$ be the tautological line bundle over $P(T^* X)$. Let $u$ be a smooth semipositive $(1,1)$-form on $X$. We say that $u$ is a \textit{smooth curvature bound} if there exists a smooth metric on $\Oo_{TX} (1)$ with curvature $c(\Oo_{TX} (1))$ such that $c(\Oo_{TX} (1))+\pi^* (u)\geq 0$. So for example if the vector bundle $TX$ is ample we can choose $u=0$.

\section{Algebraic Morse inequalities}
The special case when $F$ is nef of the theorem below was proved in~\cite{Dem-1} and the proof given here follows similar lines. See also \cite{Dem-3}.

\begin{theorem}
\label{Morse1}
Let $L$ and $F$ be holomorphic line bundles over $X$, assume that $L$ is nef and that $\codim E_{n,n}(F)\geq q+1$, and let $0\leq m\leq q$. Then we have the following algebraic Morse inequalities:
\begin{align*}
h^m(X,k(L-F))
&\leq\frac{k^n}{n!} \binom{n}{m} L^{n-m} F^m +o(k^n),\\
(-1)^m\sum_{i=0}^m (-1)^i h^i (X,k(L-F))
&\leq(-1)^m \frac{k^n}{n!}
\sum_{i=0}^m (-1)^i \binom{n}{i} L^{n-i} F^i +o(k^n).
\end{align*}
\end{theorem}

\begin{proof}
It is sufficient to prove the statement for $m=q$. Let $A$ be an ample divisor, and let $\epsilon$ be a small rational number. By replacing $L$ and $F$ with $L+\epsilon A$ and $F+\epsilon A$ respectively we may assume $L$ ample and $F$ big. Hence by~\cite{Bou}*{Proposition~3.6} $\nu^*(c_1(V),x)=\nu(\Tmin,x)$ where $\Tmin$ is a closed positive current with minimal singularities representing $c_1(F)$. Choose a positive closed $(1,1)$-form $\omega$ which is the curvature of a smooth metric on the ample line bundle $L$. The assumption that the codimension of $E_{n,n}(F)$ is at least $q+1$ is then the same thing as saying that the superlevel sets $E_c(\Tmin)$ have codimension larger or equal to $q+1$ for every positive $c$. Let us choose $u$ a smooth curvature bound. We know from~\cite{Dem-2} that there exists a sequence of closed smooth forms $T_k$ in the cohomology class of $c_1(F)$, a decreasing sequence $\lambda_k$ of positive continuous functions, and a decreasing sequence of positive real numbers $\epsilon_k$ with the following properties:
\begin{itemize}
\item the sequence $T_k$ converges weakly to $\Tmin$;
\item $\lambda_k(x)$ converges to $\nu(x,\Tmin)$ for all $x\in X$;
\item $\epsilon_k$ converges to $0$;
\item $T_k>-\lambda_k u-2\epsilon_k \omega$.
\end{itemize}

Let $c$ be a positive real number, and let $\Omega_{k,c}=\{x\in X:\lambda_k<c\}$. Let $v_k=2\epsilon_k \omega+cu$, and $w_k=\omega+v_k$. Then on $\Omega_{k,c}$ the smooth $(1,1)$-forms $T_k+v_k$ and $w_k$ are positive. Now $\omega-T_k$ is the curvature of a smooth metric on $L-F$. Let $\alpha_1\leq\dotsb\leq\alpha_n$ be the eigenvalues of $\omega-T_k$ with respect to $w_k$, so that $\alpha_j<1$ for $1\leq j\leq n$. Let $X(s)=\{x\in X:\alpha_s\leq 0,\ \alpha_{s+1}\geq 0\}$ and $X(\leq s)=\{x\in X:\alpha_{s+1}\geq 0\}$. Then on $X(s)$ we have $(-1)^s(\omega-T_k)^n\leq(-1)^s\alpha_1\dotsm\alpha_s w_k^n$. Now
\begin{equation}
\label{v-w}
\binom{n}{s} w_k^{n-s}\wedge(T_k+v_k)^s
=\binom{n}{s} w_k^{n-s}\wedge(w_k-(\omega-T_k))^s
=\sigma_s (1-\alpha)w_k^n,
\end{equation}
where $\sigma_s(1-\alpha)$ is the $s$-th elementary symmetric function in $1-\alpha_1,\dotsc,1-\alpha_n$. However, since $\alpha_j<1$ for $1\leq j\leq n$, it follows that on $X(s)$ we have
\begin{equation*}
\sigma_s(1-\alpha)\geq(1-\alpha_1)\dotsm(1-\alpha_s)
\geq(-1)^s\alpha_1\dotsm\alpha_s.
\end{equation*}
Furthermore one can easily prove by induction on $n$ \cite{Dem-1} that
\begin{equation*}
\sum_{s=0}^q (-1)^{q-s} \sigma_s (1-\alpha)
\geq\chi_{X(\leq q)} (-1)^q \alpha_1 \dotsm\alpha_n,
\end{equation*}
where $\chi_{X(\leq q)}$ is the characteristic function of the set $X(\leq q)$.

Let $\Lambda=\max_{x\in X} \nu(\Tmin,x)+1$. From the above we find
\begin{align*}
\int_{X(q)} (-1)^q (\omega-T_k)^n
&\leq\binom{n}{q} \biggl(\int_{\Omega_{k,c}} w_k^{n-q}\wedge(T_k+v_k)^q\\
&\qquad+
\int_{X\setminus\Omega_{k,c}} (\omega+\Lambda u+2\epsilon_1\omega)^{n-q}\\
&\qquad\qquad\wedge(T_k+\Lambda u+2\epsilon_1 \omega)^q \biggr),\\
\int_{X(\leq q)} (-1)^q (\omega-T_k)^n
&\leq\int_{\Omega_{k,c}}\sum_{s=0}^q (-1)^{q-s} \binom{n}{s} 
w_k^{n-s} \wedge(T_k+v_k)^s\\
&\qquad+\int_{X\setminus\Omega_{k,c}}
\sum_{s=0}^q (-1)^{q-s} \binom{n}{s}  (\omega+\Lambda u+2\epsilon_1\omega)^{n-s}\\
&\qquad\qquad\wedge(T_k+\Lambda u+2\epsilon_1\omega)^s\bigg..
\end{align*}
Consider the form $\Theta_{k,r}=(T_k+\Lambda u+2\epsilon_1\omega)^r$ for $0\leq r\leq q$. Since these forms are in a fixed cohomology class, they have bounded mass, hence up to subsequences they converge as $k\to\infty$ to an $(r,r)$-current $\Theta_r$. Since by assumption if $c>0$, then $\codim E_c(\Tmin)>q$, it follows that $\Theta_r$ has no mass on the sets $E_c$.

For $0\leq s\leq q$ the form $w_k^{n-s}\wedge(T_k+v_k)^s$ is a linear combination with smooth coefficients of the forms $(T_k+\Lambda u+2\epsilon_1\omega)^r$ for some $0\leq r\leq s$. So the forms $w_k^{n-s}\wedge(T_k+v_k)^s$ will also converge to currents with no mass on the sets $E_c(\Tmin)$. However $\bigcap_k(X\setminus\Omega_{k,c})=E_c(\Tmin)$. Therefore
\begin{equation*}
\biggl|\int_X w_k^{n-s} \wedge(T_k+v_k)^s
-\int_{\Omega_{k,c}} w_k^{n-s}\wedge(T_k+v_k)^s\biggr|\to 0,
\end{equation*}
that is to say that $\lim_{k\to\infty} \int_{\Omega_{k,c}} w_k^{n-s}\wedge(T_k+v_k)^s=(L+cu)^{n-s}(F+cu)^s$. Now we apply the standard holomorphic Morse inequalities, and conclude by letting $c\to 0$.
\end{proof}

An algebraic proof of the above theorem in the case $F$ nef was given in~\cite{A}. In case $q=1$ and $F$ nef there is a proof via the Riemann-Roch theorem due to Catanese \cite{Dem-1}.

Let $F$ be a pseudoeffective line bundle, then we have the divisorial Zariski decomposition $F=Z(F)+N(F)$ with $N(F)=\sum_{j=1}^N \nu_j D_j$. Here $\codim E_{n,n}(Z(F))\geq 2$, hence we obtain the following.

\begin{corollary}
Let $L$ and $F$ be line bundles over $X$, assume that $L$ is nef and $F$ pseudoeffective, and let $F=Z(F)+N(F)$ be the divisorial Zariski decomposition of $F$. Then
\begin{equation*}
\Vol(L-Z(F))\geq L^n-nL^{n-1} Z(F)\geq L^n-nL^{n-1} F.
\end{equation*}
\end{corollary}

\begin{proof}
Since $L$ is nef, it is sufficient to show that $\Vol(L-Z(F))\geq L^n-n L^{n-1} Z(F)$. If $Z(F)$ is a $\Q$ divisor the result follows from~Theorem~\ref{Morse1}. In the general case, we assume $L$ ample and $Z(F)$ big. Let $\epsilon$ be a small positive number, and let $\nu_{j,k}$ be an increasing sequence of rational numbers converging to $\nu_j$. Then $Z_k=Z(F)+\sum_j (\nu_j-\nu_{j,k})D_j$ is a $\Q$ divisor. Moreover it follows from the fact that the map $N$ is convex and homogeneous \cite{Bou} that $b_2(\Tmin(Z_k))<\epsilon$ for large $k$, where $\Tmin(Z_k)$ is the current with minimal singularities representing $Z_k$. Choose $u$ a smooth curvature bound. From the proof of~Theorem~\ref{Morse1} we have $\Vol(L-Z_k)\geq L^n-n(L+\epsilon u)^{n-1}(Z_k+\epsilon u)$. Now we let $\epsilon\to 0$ and $k\to\infty$.
\end{proof}

We observed in~\cite{Tra} that the corollary above is false in general if $\codim E_{n,n}(F)=1$ and we replace $Z(F)$ with $F$. However for $q=1$ there is a version of the algebraic Morse inequalities which uses the full divisorial Zariski decomposition.

\begin{theorem}
Let $L$ and $F$ be line bundles over $X$, assume that $L$ is nef and $F$ pseudoeffective, let $F=Z(F)+N(F)$ be the divisorial Zariski decomposition of $F$ with $N(F)=\sum_{j=1}^N \nu_j D_j$, and let $\{u\}$ be a nef cohomology class in $H^2(X,\R)$ such that $c_1(\Oo_{TX}(1))+\pi^*(u)$ is a nef cohomology class in $H^2(P(T^* X,\R))$. Then
\begin{equation*}
\Vol(L-F)\geq L^n-nL^{n-1} Z(F)-n\sum_{j=1}^N (L+\nu_j u)^{n-1} \nu_j D_j.
\end{equation*}
\end{theorem}

\begin{proof}
As above we see that it is sufficient to prove the theorem for $L$ ample and $F$ big. Furthermore given a positive real number $\epsilon$, by~\cite{Dem-2}*{Lemma 6.2} we can find a smooth form $u_\epsilon$ in the class $\{u\}$ such that $u_\epsilon+\epsilon A$ is a smooth curvature bound. Here $A$ is an ample divisor. If we prove the theorem with $\{u\}$ replaced by $\{u\}+\epsilon A$ then it is sufficient to let $\epsilon\to 0$ to obtain the general result. So we can further assume that $u$ is a smooth curvature bound. Again by~\cite{Bou}*{Proposition~3.6} if $\Tmin$ is a positive current with minimal singularities representing $F$ then $\nu^*(c_1(F),x)=\nu(\Tmin,x)$, in particular the divisorial Zariski decomposition $F=Z(F)+\sum_{j=1}^N \nu_j D_j$ is just the Siu decomposition of the current $\Tmin$ \cite{Siu-2}. For $1\leq j\leq N$, let us choose an increasing sequence of rational numbers $\nu_{j,k}$ converging to $\nu_j$ as $k\to\infty$. Let $Z_k=Z(F)+\sum_j (\nu_j -\nu_{j,k}) D_j$ and $N_k=\sum_j\nu_{j,k}D_j$, then $Z_k$ and $N_k$ are effective $\Q$ divisors such that $F=Z_k+N_k$. Moreover $\lim_{k\to\infty} b_2 (\Tmin-\sum_j \nu_{j,k} [D_j])=0$. Let us for the moment fix a positive integer $k$. By~\cite{Dem-4}*{proof of~Proposition~9.1} there is a sequence $\Delta_s$ of effective divisors and a sequence $m_s$ of positive integers such that as $s\to\infty$ we have that:
\begin{itemize}
\item the integration currents $\frac{1}{m_s} [\Delta_s]$ converge to $\Tmin-\sum_j \nu_{j,k} [D_j]$;
\item $b_2(\frac{1}{m_s} [\Delta_s])$ converges to $b_2(\Tmin-\sum\nu_{j,k} [D_j])$;
\item for $1\leq j\leq N$ the generic Lelong number $\alpha_{j,s}$ of $\frac{1}{m_s} [\Delta_s]$ along $D_j$ converges to the generic Lelong number $\nu_j-\nu_{j,k}$ of $\Tmin-\sum_i\nu_{i,k}[D_i]$, along $D_j$.
\end{itemize}

So we can write $\frac{1}{m_s} \Delta_s=\frac{1}{m_s} \Delta'_s +\sum_j \alpha_{j,s} D_j$ where $\codim(\Delta'_s \cap D_j)\geq 2$ for $1\leq j\leq N$. Then $\limsup_{s\to\infty} b_2(\frac{1}{m_s} [\Delta'_s])\leq b_2 (\Tmin-\sum_j \nu_{j,k} [D_j])$. Now we fix $s$ and replace the divisor $F$ by the divisor $F'=\frac{1}{m_s} \Delta'_s +\sum_j (\nu_{j,k}+\alpha_{j,s}) D_j$. Let $U_h (0)$ and $U_h (j)$, for $1\leq j\leq N$, be decreasing fundamental systems of open neighborhoods of $\Delta'_s$ and $D_j$ respectively. Let $\Delta_{0,h}$ and ${D_{j,h}}$ sequences of smooth forms approximating for $h\to\infty$ the currents $\frac{1}{m_s}[\Delta'_s]$ and $(\nu_{j,k}+\alpha_{j,s})[D_j]$ as in the proof of~Theorem~\ref{Morse1} above. By~\cite{Tra}*{Lemma~2.4} for large $h$ the forms $\Delta_{0,h}$ and $D_{j,h}$ can be taken to have support in $U_h (0)$ and in $U_h (j)$ respectively. Let $\Theta_h=\Delta_{0,h}+\sum_j D_{j,h}$, then $\omega-\Theta_h$ is the curvature of a smooth metric on $L- F'$. Observe that if we look at such a metric for large $h$, then $X_1\subseteq U_h (0)\cup\bigcup_j U_h (j)$ and for $K\subseteq X\setminus\Delta'_s \setminus\bigcup_j D_j$ compact we may assume that $K\cap(U_h (0)\cup\bigcup_j U_h (j))=\emptyset$. Therefore we have $\int_K(\omega-\Theta_h)^n \geq\int_K\omega^n$. It follows that
\begin{equation*}
\liminf_{h\to\infty} \int_{X(0)} (\omega-\Theta_h)^n \geq L^n.
\end{equation*}
Now we use the notation of~Theorem~\ref{Morse1}. Let us fix a positive real number $c$ and look at~formula~\eqref{v-w} with $s=1$ on the open set  $X_1\cap\Omega_{h,c}$. We find that
\begin{equation*}
nw_h^{n-1}\wedge(\Theta_h+v_h)=\sum_{j=1}^n (1-\alpha_j)w_h^n
\geq-\alpha_1 w_h^n \geq-(\omega-\Theta_h)^n.
\end{equation*}
Moreover $w_h^{n-1} \wedge(\Theta_h+v_h)$ converges weakly to $(\omega+cu)^{n-1} \wedge(\Tmin+cu)$ as $h\to\infty$. For $h$ increasing the measure of $X(1)$ tends to $0$, therefore
\begin{equation*}
\lim_{h\to\infty} \int_{X(1)} (\omega+cu)^{n-1} \wedge(cu)=0.
\end{equation*}
On the other hand the intersection of any two distinct divisors $D_{j_1}$ and $D_{j_2}$ and the intersection of each $D_j$ with $\Delta'_s$ has codimension at least $2$. Now for fixed $j$, with $0\leq j\leq N$, we let $V_h (j)=\bigcup_{i\neq j} (U_h (i)\cap U_h (j))$ for $0\leq i\leq N$. Then as in the proof of~Theorem~\ref{Morse1} we see that
\begin{equation*}
\lim_{h\to\infty} \int_{X(1)\cap V(j)} -(\omega-\Theta_h)^n =0
\end{equation*}
and that if we take $c_0 >b_2 (\frac{1}{m_s} \Delta'_s)$ then
\begin{equation*}
\begin{split}
\limsup_{h\to\infty}
\int_{X(1)\cap U_h (0)\setminus V_h (0)} -(\omega-\Theta_h)^n
&\leq n(L+c_0 u)^{n-1} \frac{1}{m_s} [\Delta'_s]\\
&\leq n(L+c_0 u)^{n-1} \frac{1}{m_s} [\Delta_s].
\end{split}
\end{equation*}
Similarly if we take $c_j>\nu_{j,k}+\alpha_{j,s}$ for $1\leq j\leq N$ we obtain
\begin{equation*}
\limsup_{h\to\infty}
\int_{X(1)\cap U_h (j)\setminus V_h (j)} -(\omega-\Theta_h)^n
\leq n(L+c_j u)^{n-1} (\nu_{j,k} +\alpha_{j,s})[D_j].
\end{equation*}
Putting everything together and using the standard holomorphic Morse inequalities, we obtain:
\begin{multline*}
h^0 (X,k(L-F'))-h^1 (X,k(L-F'))\\
\begin{aligned}
&\geq\frac{k^n}{n!}\biggl(L^n-n(L+c_0 u)^{n-1} \frac{1}{m_s} [\Delta_s]\\
&\qquad
-n\sum_{j=1}^N (L+c_j u)^{n-1} (\nu_{j,k}+\alpha_{j,s})[D_j]\biggr)+o(k^n).
\end{aligned}
\end{multline*}
We now let $c_j$ converge to $\nu_{j,k}+\alpha_{j,s}$ and $c_0$ to $b_2(\frac{1}{m_s} [\Delta'_s])$. To conclude, first we let $s\to\infty$, then $k\to\infty$ and use the continuity of the volume.
\end{proof}

\end{document}